\documentclass[centertags,12pt]{amsart}
\usepackage{latexsym}
\usepackage{amssymb}

\numberwithin{equation}{section}
\makeatletter
\renewcommand{\subsection}{\@startsection
{subsection}{2}{0mm}{\baselineskip}{-0.25cm}
{\normalfont\normalsize\bf}}
\makeatother

\newtheorem{theorem}{Theorem}[section]
\newtheorem{proposition}[theorem]{Proposition}
\newtheorem{corollary}[theorem]{Corollary}
\newtheorem{lemma}[theorem]{Lemma}
\theoremstyle{definition}
\newtheorem{remark}[theorem]{Remark}
\newtheorem{definition}[theorem]{Definition}
\newtheorem{problem}[theorem]{Problem}

\def\P{\mathbf P}

\def\cA{\mathcal A}

\def\cC{\mathcal C}
\def\cD{\mathcal D}

\def\cK{\mathcal K}
\def\cX{\mathcal X}
\def\cY{\mathcal Y}

\def\sq{\sqrt{q}}
\def\fq{\mathbf F_q}

\def\supp{{\rm Supp}}

\def\dim{{\rm dim}}
\def\deg{{\rm deg}}
\def\det{{\rm det}}
\def\frx{{\mathbf Fr}_{\mathcal X}}

\begin{document}

\author[Giulietti]{M. Giulietti}
\author[Pambianco]{F. Pambianco}
\author[Torres]{F. Torres}
\author[Ughi]{E. Ughi}\thanks{1991 Math. Subj. Class.: 
Primary 05B, Secondary 14G}

\title[Large complete arcs]{On large complete arcs: odd case}
\address{Dipartimento di Matematica - Universit\`a degli Studi di Perugia,
Via Vanvitelli 1 - 06123 Perugia (Italy)}
\email{giuliet@dipmat.unipg.it}
\email{fernanda@dipmat.unipg.it}
\address{IMECC-UNICAMP, Cx. P. 6065, Campinas, 13083-970-SP, Brazil}
\email{ftorres@ime.unicamp.br}
\address{Dipartimento di Matematica - Universit\`a degli Studi di Perugia,
Via Vanvitelli 1 - 06123 Perugia (Italy)}
\email{ughi@dipmat.unipg.it}

     \begin{abstract} 
An approach for the computations of upper bounds on the size of large
complete arcs is presented. We obtain in particular geometrical
properties of irreducible envelopes associated to a second largest 
complete arc provided that the order of the underlying field is 
large enough.
    \end{abstract}

\maketitle

   \section{Introduction}\label{s1}
A {\em $k$-arc} in the projective plane $\P^2(\fq)$, where $\fq$ is the
finite field with $q$ elements, is a set of $k$ points no three of which
are collinear. An arc is {\em complete} if it is not properly contained in
another arc. For a given $q$, a basic problem in Finite Geometry is to
find the values of $k$ for which a complete $k$-arc exists. For a $k$-arc  
$\cK$ in $\P^2(\fq)$, Bose \cite{bose} showed that
$$
k\le m_2(q):=\begin{cases}
q+1 & \text{if $q$ is odd}\, ,\\
q+2 & \text{otherwise}\, .
\end{cases}
$$  
For $q$ odd the bound $m_2(q)$ is attained if and only if $\cK$ is an
irreducible conic \cite{segre}, \cite[Thm. 8.2.4]{h}. For $q$ even the
bound is attained by the union of an irreducible conic and its nucleus,
and not every $(q+2)$-arc arises in this way; see \cite[\S8.4]{h}. Let
$m_2'(q)$ denote the second largest size that a complete arc in
$\P^2(\fq)$ can have. Segre \cite{segre}, \cite[\S10.4]{h} showed that 
    \begin{equation}\label{eq1.1}
m_2'(q)\le\begin{cases}
q-\frac{1}{4}\sq+\frac{7}{4} & \text{if $q$ is odd},\\
q-\sq+1                      & \text{otherwise}.
\end{cases}
    \end{equation}
Besides small $q$, namely $q\le 29$ \cite{chao-kaneta}, \cite{h},
\cite{h-storme}, the only case where $m_2'(q)$ has been determinated is
for $q$ an even square.  Indeed, for $q$ square, examples of complete
$(q-\sq+1)$-arcs \cite{boros}, \cite{coss}, \cite{ebert},
\cite{fisher-h-thas}, \cite{kestenband} show that
  \begin{equation}\label{eq1.2}
m_2'(q)\ge q-\sq+1\, ,
   \end{equation}
and so the bound (\ref{eq1.1}) for $q$ an even square is sharp. This
result has been recently extended by Hirschfeld and Korchm\'aros
\cite{hk3} who showed that the third largest size that a complete arc
can have is bounded from above by $q-2\sq+6$. 

If $q$ is not a square, Segre's bounds were notably improved by
Voloch \cite{v1}, \cite{v2} (see \S3 here). 

If $q$ is odd, Segre's bound was slightly improved to $m_2'(q)\le
q-\sq/4+25/16$ by Thas \cite{thas}. If $q$ is an odd square and large
enough, Hirschfeld and Korchm\'aros \cite{hk2} significantly improved the
bound to
   \begin{equation}\label{eq1.3}
m_2'(q)\le q-\frac{1}{2}\sq+\frac{5}{2}\, . 
   \end{equation}
The two last bounds suggest the following problem, which seems to be
difficult and has remained open since the 60's.
   \begin{problem}\label{prob1.1}
For $q$ an odd square, is it true that $m_2'(q)=q-\sq+1$?
   \end{problem}
The answer is no for $q=9$ and yes for $q=25$ \cite{chao-kaneta},
\cite{h}, \cite{h-storme}. So Problem \ref{prob1.1} is indeed open for 
$q\ge 49$.

In this paper we investigate irreducible components of the envelope
associated to large arcs in $\P^2(\fq)$. Such components will be called
{\em irreducible envelopes} and their existence is related to the
existence of certain rational points which will be called {\em special
points}, see \S2. This set up allows us to prove a general bound for the
size of a complete arc (Proposition \ref{prop3.1}) which depends on $q$
and the 4th positive $\fq$-Frobenius order of the linear series obtained
from quadrics in $\P^2(\bar\fq)$ defined on any irreducible envelope. 
From this result, for $q$ odd and not a square, we recover the bounds on
the size of arcs that were established so far in the literature (Lemma
\ref{lemma3.1}, Lemma \ref{lemma3.5}). For $q$ an odd square the best that
our approach gives is another proof of Segre's bound. 

Our research was inspired and motivated by the papers of Voloch
\cite[\S4]{v2} and Hirschfeld--Korchm\'aros \cite{hk1}, \cite{hk2}. In
fact our results are implicitly contained in such works and this paper can
be considered as a set of footnotes to those. Nonetheless, the main
contribution of this paper are the following. 
   \begin{enumerate}
\item[(I)] We explicitly determinate the type of curves (see (III) below) 
associated to complete large arcs (Proposition \ref{prop4.1} and
Proposition \ref{prop4.2})  whenever the underlying field is large enough; 

\item[(II)] We give a systematic account of how to bound the size of
complete arcs by means of St\"ohr-Voloch's approach to the Hasse-Weil
bound \cite{sv}; 

\item[(III)] We provide motivation for the study of irreducible plane
curves over $\fq$ whose $\fq$ non-singular model is classical for the
linear series $\Sigma_1$ obtained from lines and whose
$\Sigma_2:=2\Sigma_1$ orders 
are $0,1,2,3,4,\epsilon_5$  and whose 
 $\fq$-Frobenius orders for $\Sigma_2$ are $0,1,2,3,\nu_4$, where
$\epsilon_5=\nu_4\in \{\sqrt{q/p},\sq, \sq/3,3\sq\}$. See \S4 here.
     \end{enumerate}

Finally, for the convenience of the reader, we include an appendix
containing basic facts from Weierstrass points and Frobenius orders based
on St\"ohr-Voloch's paper \cite{sv}.
   \section{Special points and irreducible envelopes}\label{s2}
Throughout this section $\cK$ will be an arc in $\P^2(\fq)$. Segre
associates to $\cK$ a plane curve $\cC$ in the dual plane of
$\P^2(\bar\fq)$, where $\bar\fq$ denotes the algebraic closure of $\fq$.
This curve is defined over $\fq$ and it is called {\em
the envelope of} $\cK$. For $P\in \P^2(\bar\fq)$, let $\ell_P$ denote the
corresponding line in the dual plane. The following result summarize the
main properties of $\cC$ for the odd case.
    \begin{theorem}\label{thm2.1} If $q$ is odd, then the following
statements hold:
    \begin{enumerate}

\item The degree of $\cC$ is $2t$, with $t=q-k+2$ being the number of
1-secants through a point of $\cK$.

\item All $kt$ of the  $1$-secants of $\cK$  belong to $\cC$.

\item Each $1$-secant $\ell$ of $\cK$ through a point $P\in\cK$ is counted
twice in the intersection of $\cC$ with $\ell_P$, i.e. 
$I(\ell, \cC\cap \ell_P)=2$.

\item The curve $\cC$ contains no $2$-secant of $\cK$.

\item The irreducible components of $\cC$ have multiplicity at most two,
and $\cC$ has at least one component of multiplicity one.

\item The arc $\cK$ is incomplete if and only if $\cC$ admits a linear
component over $\fq$. The arc $\cK$ is a conic if and only it is complete
and $\cC$ admits a quadratic component over $\fq$. 

\end{enumerate}
   \end{theorem}
   \begin{proof} See \cite{segre}, \cite[\S10]{h}.
    \end{proof}
We recall that a non-singular point $P$ of a plane curve $\cA$ is 
called an {\em inflexion point of $\cA$} if $I(P, \cA\cap\ell)>2$, with
$\ell$ being the tangent line of $\cA$ at $P$. We introduce the following
terminology:
  \begin{definition}\label{def2.1} A point $P_0$ of $\cC$ is called {\em
special} if the following conditions hold:
   \begin{enumerate} 
\item[(i)] it is non-singular; 
\item[(ii)] it is $\fq$-rational;
\item[(iii)] it is not an inflexion point of $\cC$.
   \end{enumerate}
Then, by (i), a special point $P_0$ belongs to an 
unique irreducible component of the envelope which will be called {\em the
irreducible envelope} associated to $P_0$ or {\em an irreducible envelope
of $\cK$}. 
   \end{definition}
\begin{lemma}\label{lemma2.1}  Let $\cC_1$ be an irreducible envelope of
$\cK$. Then
\begin{enumerate}
\item $\cC_1$ is defined over $\fq;$
\item  if $q$ is odd and the arc is not a conic 
and complete, then the degree of $\cC_1$ is at least three.
\end{enumerate}
\end{lemma}
   \begin{proof} (1) Let $\cC_1$ be associated to $P_0$, let $\mathbf \Phi$ be
the Frobenius morphism (relative to $\fq$) on the dual plane of
$\P^2(\bar\fq)$, and suppose that $\cC_1$ is not defined over $\fq$. Then,
since the envelope is defined over $\fq$ and $P_0$ is $\fq$-rational,
$P_0$ would belong to two different components of the envelope, namely
$\cC_1$ and ${\mathbf \Phi}(\cC_1)$. This is a contradiction because the
point is non-singular. 

(2) This follows from Theorem \ref{thm2.1}(6). 
\end{proof}
The next result will show that special points do exist provided that $q$ is
odd and the arc is large enough.
     \begin{proposition}\label{prop2.1} Let $\cK$ be an arc in $\P^2(\fq)$  
of size $k$ such that $k>(2q+4)/3$. If $q$ is odd, then the envelope 
$\cC$ of $\cK$ has special points.
     \end{proposition}
  \begin{remark}\label{rem2.1} The hypothesis $k>(2q+4)/3$ in the
proposition is equivalent
to $k>2t$, with $t=q-k+2$. Also, under this hypothesis, the envelope
$\cC$ is uniquely determined by $\cK$, see \cite[Thm. 10.4.1(i)]{h}.
  \end{remark}
To prove Proposition \ref{prop2.1} we need the following lemma, for which
we could not find a reference. 
      \begin{lemma}\label{lemma2.2} 
Let $\cA$ be a plane curve defined over $\bar\fq$ and suppose that it 
has no multiple components. Let $\alpha$ be the degree of $\cA$ and $s$
the number of its singular points. Then,
$$
s\le \binom{\alpha}{2}\, ,
$$
and equality holds if $\cA$ consists of $\alpha$ lines no three 
concurrent.  
      \end{lemma}
      \begin{proof}
That a set of $\alpha$ lines no three concurrent satisfies the bound is
trivial. Let $G=0$ be the equation of $\cA$, let $G=G_1\ldots G_r$ be the 
factorization of $G$ in $\bar\fq[X,Y]$, and let $\cA_i$ be the curve given by
$G_i=0$. For 
simplicity we assume $\alpha$ even, say $\alpha=2M$. Setting
$\alpha_i:=\deg(G_i)$, 
$i=1, \ldots, r$ and $I:=\sum_{i=1}^{r-1}\alpha_i$ we have
$\alpha_r=2M-I$. The
singular points of $\cA$ arise from the singular points of each component
or from the points in $\cA_i \cap \cA_j $, $i \neq j$. Recall that an
irreducible plane curve of degree $d$ has at most $\binom{d-1}{2}$ singular
points, and that $\#\cA_i \cap \cA_j \le a_ia_j$, $i \neq j$ (B\'ezout's
Theorem). So
\begin{equation*}
\begin{split}
s  & \le  \sum_{i=1}^{r-1}\binom{\alpha_i
-1}{2}+\binom{2M-I-1}{2}+\sum_{1\le
i_1<i_2\le r-1}\alpha_{i_1}\alpha_{i_2}+\sum_{i=1}^{r-1}(2M-I)\alpha_i\\
 & =
\sum_{i=1}^{r-1}\frac{\alpha_i^2-3\alpha_i+2}{2}+\frac{4M^2-4MI+I^2-6M+3I+2}{2}
+\sum_{1\le i_1 < i_2 \le r-1}\alpha_{i_{1}}\alpha_{i_{2}}+(2M-I)I\\
 & =
\frac{1}{2}[
\sum_{i=1}^{r-1}\alpha_{i}^{2}-3I+2(r-1)+4M^2-4MI+I^2-6M+3I+2+\\
{} &\quad \ 2\sum_{1\le i_1 
< i_2 \le r-1}\alpha_{i_{1}}\alpha_{i_{2}}+4MI-2I^2]\\
 & \le 2 M^2 - 3M + \alpha = 2M^2 -M\, .
\end{split}
\end{equation*}
          \end{proof}
          \begin{proof} ({\em Proposition \ref{prop2.1}}) Let $F=0$ be the
equation of $\cC$ over $\fq$. By Theorem \ref{thm2.1}(5), $F$ admits a 
factorization in $\bar\fq[X,Y,Z]$ of type
$$
G_1\ldots G_r H_1^2\ldots H_s^2\, ,
$$
with $r\ge 1$ and $s\ge 0$. Let $\cA$ be the plane curve given by 
  $$
G:=G_1\ldots G_r=0\, .
  $$
Then $\cA$ satisfies the hypothesis of Lemma \ref{lemma2.1} and it has  
even degree by Theorem \ref{thm2.1}(1). From Theorem \ref{thm2.1}(3)  
and B\'ezout's theorem,  for each line $\ell_P$ (in the dual plane)  
corresponding to a point $P\in \cK$, we have  
$$
\#(\cA \cap \ell_P) \ge M\, ,
$$
where $2M = \deg(G)$, and so at least $kM$ points corresponding 
to unisecants of $\cK$ belong to $\cA$. Since $k>2t$ (see
Remark \ref{rem2.1}) and $2t\ge 2M$, then $kM>2M^2$ and from Lemma 
\ref{lemma2.1} we have that at least one
of the unisecant points in $\cA$, says $P_0$, is non-singular. Suppose
that $P_0$ goes through $P\in \cK$. The point $P_0$ is clearly
$\fq$-rational and 
$P_0$ is not a point of the curve of equation $H=0$: 
otherwise $I(P_0, \cC\cap\ell_P) > 2$ (see Theorem \ref{thm2.1}(3)). Then,  
$I(P_0, \cC\cap\ell_P)=I(P_0,\cA\cap\ell_P)=2$ and so $\ell_P$ is
the tangent of $\cC$ at $P_0$. Therefore $P_0$ is not an inflexion 
point of $\cC$, and the proof of Proposition \ref{prop2.1} is complete. 
             \end{proof}
Let $\cC_1$ be an irreducible envelope associated to a special point
$P_0$, and 
$$
\pi: \cX\to \cC_1\, ,
$$
the normalization of $\cC_1$. Then by Lemma \ref{lemma2.1}(1) we can
assume that $\cX$ and $\pi$ are defined over $\fq$. In particular, the
linear series $\Sigma_1$ on $\cX$ obtained by the pullback of lines of
$\P^2(\bar\fq)^*$, 
the dual of  
$\P^2(\bar\fq)$, is $\fq$-rational. Also, there is just one point
$\tilde P_0\in \cX$ such that $\pi(\tilde P_0)=P_0$. For basic facts on
orders and Frobenius orders the reader is referred to \cite{sv} or the
appendix here.
   \begin{lemma}\label{lemma2.3} Let $q$ be odd. Then,
\begin{enumerate}
  \item the $(\Sigma_1,\tilde P_0)$-orders are $0, 1, 2;$
  \item the curve $\cX$ is classical with respect to $\Sigma_1$.
\end{enumerate}
    \end{lemma}
\begin{proof} (1) This follows from the proof of Proposition \ref{prop2.1}.

(2) This follows from Item (1) and (W1) in the appendix.
\end{proof}
   \begin{remark}\label{rem2.2} The hypothesis $q$ odd in Lemma
\ref{lemma2.3} (as well as in Proposition \ref{prop2.1}) 
is necessary. In fact, from \cite{fisher-h-thas} and \cite{thas} follow
that the envelope associated to the cyclic $(q-\sq+1)$-arc, with $q$ an even
square, is irreducible and $\fq$-isomorphic to the plane curve
$XY^{\sq}+X^{\sq}Z+YZ^{\sq}=0$ which is not $\Sigma_1$-classical. 
   \end{remark}
Next consider the following sets:
\begin{align*}
\cX_1(\fq):= & \{P\in \cX: \pi(P)\in \cC_1(\fq)\}\, ,\\
\cX_{11}(\fq):= & \{P\in \cX_1(\fq): j^1_2(P)=2j^1_1(P)\}\, ,\\
\cX_{12}(\fq):= & \{P\in \cX_1(\fq): j^1_2(P)\neq 2j^1_1(P)\}\, ,
\end{align*}
and the following numbers:
   \begin{equation}\label{eq2.1}
M_q=M_q(\cC_1):=\sum_{P\in \cX_{11}(\fq)} j^1_1(P)\, , \qquad
M_q'=M_q'(\cC_1):=\sum_{P\in \cX_{12}(\fq)} j^1_1(P)\, ,
   \end{equation}
where $0<j^1_1(P)<j^1_2(P)$ denotes the $(\Sigma_1,P)$-order sequence. We
have that 
$$
M_q+M_q'\ge \#\cX_1(\fq) \ge \#\cX(\fq)\qquad\text{and}\qquad
\#\cX_1(\fq)\ge \#\cC_1(\fq)\, .
$$
    \begin{proposition}\label{prop2.2} Let $\cK$ be an arc of size $k$ and  
$d$ the degree of an irreducible envelope of 
$\cK$. For $M_q$ and $M_q'$ as above we have
$$
2M_q+M_q'\ge kd\, .
$$
    \end{proposition}
To prove the proposition we first prove the following
	\begin{lemma}\label{lemma2.4}
Let $\cK$ be an arc and $\cC_1$ an irreducible envelope of $\cK$.  Let $Q
\in \cK$ and $\cA_Q$ be the set of points of $\cC_1$ corresponding to
unisecants of $\cK$ passing through $Q$. Let $u:=\#\cA_Q$ and $v$ be the
number of points in $\cA_Q$ which are non-singular and inflexion points
of $\cC_1$. Then
$$
2(u-v)+v \ge d\, ,
$$
where $d$ is the degree of $\cC_1$.
	\end{lemma}
       \begin{proof} 
Let $P'\in \cA_Q$.  Suppose that it is  non-singular and an inflexion point
of $\cC_1$. Then, from Theorem \ref{thm2.1}(3) and the definition of
$\cA_Q$, we have that $\ell_Q$ is not the tangent line of $\cC_1$ at $P'$,
i.e. we have that $I(P', \cC_1\cap \ell_Q)=1$. Now suppose that $P'$ is
either singular or a non-inflexion point of $\cC_1$. Then from Theorem
\ref{thm2.1}(3) we have $I(P',\cC_1\cap \ell_Q)\le 2$ and the result
follows from B\'ezout's theorem applied to $\cC_1$ and $\ell_Q$. 
       \end{proof}
{\em Proof of Proposition \ref{prop2.2}.} For $Q\in \cK$ let $\cA_Q$ be
as in Lemma \ref{lemma2.4} and set 
$$
\cY_Q:=\{P \in \cX_1(\fq) :  \pi(P) \in \cA_Q \}\, .
$$
We claim that
$$
m(Q):=2\sum_{P \in \cX_{11}(\fq) \cap \cY_Q}j^1_1(P)+
\sum_{P \in \cX_{12}(\fq) \cap \cY_Q}j^1_1(P) \ge d\, .
$$
This claim implies the proposition since, from Theorem \ref{thm2.1}(4),
$$
\cY_Q \cap \cY_{Q_1} = \emptyset\qquad 
\text{whenever}\qquad Q \neq Q_1\, .
$$
To prove the claim we distinguish four types of points in $\cY_Q$, namely
\begin{align*}
\cY_Q^1:= & \{P\in \cY_Q: \text{$\pi(P)$ is non-singular and non-
inflexion point of $\cC_1$}\}\, , \\
\cY_Q^2:= & \{P\in \cY_Q:  \text{$\pi(P)$ is a non-singular 
inflexion point of $\cC_1$}\}\, ,\\
\cY_Q^3:= & \{P\in \cY_Q: \text{$\pi(P)$ is a singular
point of $\cC_1$ such that $\# \pi^{-1}(\pi(P))=1$}\}\, ,\\
\cY_Q^4:= & \{P\in \cY_Q: \text{$\pi(P)$ is a singular 
point of $\cC_1$ such that $\#\pi^{-1}(\pi(P)) >1$}\}\, .
\end{align*}
Observe that $\cY_Q^1 \subseteq \cX_{11}(\fq)$ and so
$$
m(Q)\ge 2\sum_{P\in \cY_Q^1}j^1_1(P)+
\sum_{P\in \cY_Q^2}j^1_1(P)+\sum_{P\in \cY_Q^3}j^1_1(P)+\sum_{P\in
\cY_Q^4}j^1_1(P)\, .
$$
Since $j^1_1(P) > 1$ for all $P\in \cY_Q^4$, the above inequality becomes
$$
m(Q)\ge 2\#\cY_Q^1+2\#\cY_Q^4+\#\cY_Q^3+\#\cY_Q^2\, .
$$ 
Therefore, as to each singular non-cuspidal point of
$\cC_1$ in $\cA_Q$ corresponds at least two points in $\cY_Q^3$, it
follows that
\begin{align*}
m(Q) & \ge 2\#\{P' \in  \cA_Q : \text{$P'$ is either singular or not 
an inflexion point of $\cC_1$} \}+\\ 
{} &\quad  \#\{P' \in  \cA_Q : \text{$P'$ is a nonsingular inflexion
point of $\cC_1$}\}\, .
\end{align*}
Then the claim follows from Lemma \ref{lemma2.4} and the proof of
Proposition \ref{prop2.2} is complete. 
   \section{Bounding the size of an arc}\label{s3} 
Throughout the whole section we fix the following notation:
    \begin{itemize}
\item $q$ is a power of an odd prime $p$;

\item $\cK$ is a complete arc of size $k$ such that $(2q+4)/3<k\le
m_2'(q)$; therefore the degree of any irreducible envelope of $\cK$ has
at least degree three;

\item $P_0$ is an special point of the envelope $\cC$ of $\cK$
and the plane curve $\cC_1$ of degree $d$ is an irreducible envelope
associated to $P_0$; 

\item $\pi: \cX\to \cC_1$ is the normalization of $\cC_1$ which is
defined over $\fq$; as a matter of terminology, $\cX$ will be also called
an irreducible envelope of $\cK$.

\item $\tilde P_0$ is the only point in $\cX$ such that
$\pi(\tilde P_0)=P_0$; $g$ is the genus of $\cX$ (so that $g\le
(d-1)(d-2)/2$);

\item The symbols $\cX_1(\fq)$, $M_q$ and $M_q'$ are as in \S2; 

\item $\Sigma_1$ is the linear series $g^2_d$ on $\cX$ obtained
from the pullback of lines of $\P^2(\bar\fq)^*$; $\Sigma_2$ is the linear
series $g^5_{2d}$ on $\cX$ obtained from the
pullback of conics of $\P^2(\bar\fq)^*$, i.e. $\Sigma_2=2\Sigma_1$
(notice that $\dim(\Sigma_2)=5$ because $d\ge 3$);

\item $S$ is the $\fq$-Frobenius divisor associated to $\Sigma_2$;

\item $j_5(\tilde P_0)$  is the 5th positive 
$(\Sigma_2,\tilde P_0)$-order;  $\epsilon_5$ is the 5th positive 
 $\Sigma_2$-order; 
$\nu_4$ is the 4th positive $\fq$-Frobenius order of $\Sigma_2$. 
    \end{itemize}

We apply the appendix to both $\Sigma_1$ and $\Sigma_2$. We have already
noticed that the $(\Sigma_1,\tilde P_0)$-orders, as well as the
$\Sigma_1$-orders, are 0,1 and 2; see Lemma \ref{lemma2.3}. 
Then, the $(\Sigma_2, \tilde P_0)$-orders are 
0,1,2,3,4 and $j_5(\tilde P_0)$, with $5\le j_5(\tilde P_0)\le 2d$, and
the $\Sigma_2$-orders are 0,1,2,3,4 and $\epsilon_5$ with $5\le \epsilon_5\le 
j_5(\tilde P_0)$; cf. \cite[p. 464]{garcia-voloch}.

Then, we compute the $\fq$-Frobenius orders of $\Sigma_2$. We apply (F3)
in the appendix to $\tilde P_0$ and conclude that this sequence is 0,1,2,3
and $\nu_4$, with
$$
\nu_4\in\{4,\epsilon_5\}\, .
$$
Therefore (see appendix)
    \begin{align*}
\deg(S) & =(6+\nu_4)(2g-2)+(q+5)2d\, ,\\
\intertext{and}
v_P(S)& \ge 5j^2_1(P), \qquad \text{for each $P\in \cX_1(\fq)$}\, ,
    \end{align*}
where $j^2_1(P)$ stands for the first positive $(\Sigma_2,P)$-order. Since
$j^2_1(P)$ is equal to the first positive $(\Sigma_1,P)$-order (cf.
\cite[p. 464]{garcia-voloch}), we then have 
$$
\deg(S)\ge 5(M_q+M_q')\, ,
$$
where $M_q$ and $M_q'$ were defined in (\ref{eq2.1}). Then,  
taking into consideration the following facts:
   \begin{enumerate}
\item $2M_q+M_q'\ge kd$ (Proposition \ref{prop2.2}),
\item $2g-2\le d(d-3)$,
\item $\nu_4\le j_5(\tilde P_0)-1\le 2d-1$ ((F3) appendix), and
\item $d\le 2t=2(q+2-k)$ (Theorem \ref{thm2.1}(1)),
    \end{enumerate}
we obtain the following. 
   \begin{proposition}\label{prop3.1}
Let $\cK$ be a complete arc of size $k$ such that $(2q+4)/3< k\le
m_2'(q)$. Then
$$
k\le \min\{q-\frac{1}{4}\nu_4+\frac{7}{4}\, ,\   
\frac{28+4\nu_4}{29+4\nu_4}q+\frac{32+2\nu_4}{29+4\nu_4}\}\, ,
$$
where $\nu_4$ is the 4th positive $\fq$-Frobenius order of the linear 
series $\Sigma_2$ defined on an irreducible envelope of $\cK$.
    \end{proposition}
Now consider separately the cases $\nu_4=4$ and $\nu_4=\epsilon_5$.

{\bf 1.  $\nu_4=4$.} 

In this case, the corresponding irreducible
envelope will be called {\em Frobenius 
classical}. Proposition \ref{prop3.1} becomes the following.
    \begin{lemma}\label{lemma3.1}  
Let $\cK$ be a complete arc of size $k$ such that $(2q+4)/3<k\le m_2'(q)$.   
Suppose that $\cK$ admits a Frobenius classical
irreducible envelope. Then 
$$
k\le \frac{44}{45}q+\frac{40}{45}\, .
$$     
    \end{lemma}
This lemma holds in the following cases:

(3.1.1) Whenever $q=p$ is an odd prime: Voloch's bound \cite{v2}; 

(3.1.2) The arc is cyclic of Singer type whose size $k$ satisfies  
$2k\not\equiv -2,1,2,4 \pmod{p}$, where $p>5$; see Giulietti's paper
\cite{giu}.
 
For the sake of completeness we prove (3.1.1)
   \begin{proof} ({\em Item (3.1.1)}) Let 
$\cC_1$ be an irreducible envelope of
$\cK$ and $d$ the degree of $\cC_1$. If $p< 2d$, then $p<4t=4(p+2-k)$
so that $k<(3p+8)/4$ and the result follows. So let $p\ge 2d$. Then from
\cite[Corollary 2.7]{sv} we have that $\cC_1$ is Frobenius classical and
(3.1.1) follows from Proposition \ref{prop3.1}. 
    \end{proof}
Next we show that, for $q$ square and $k=m_2'(q)$, Lemma \ref{lemma3.1} 
is possible only for $q$ small.
    \begin{corollary}\label{cor3.1} Let $\cK$ be an arc of size
$m_2'(q)$ and suppose that $q$ is a square. Then,
   \begin{enumerate}
\item if $q>9$, $\cK$ has irreducible envelopes;
\item if $q>43^{2}$, any irreducible envelope of $\cK$ is Frobenius
non-classical. 
   \end{enumerate}
    \end{corollary}
    \begin{proof} (1) As we mentioned in (\ref{eq1.2}), $m_2'(q)\ge
q-\sq+1$. Since  $q-\sq+1>(2q+4)/3$ for $q>9$, Item (1) follows from
Proposition \ref{prop2.1}.

(2) If would exist a Frobenius classical irreducible envelope of $\cK$,
then from Lemma \ref{lemma3.1} and (\ref{eq1.2}) we would have
$$
q-\sq+1\le m_2'(q)\le 44q/45+40/45\, .
$$
so that $q\le 43^2$. 
   \end{proof}
{\bf 2. $\nu_4=\epsilon_5$.} 

Here, from \cite[Corollary 3]{garcia-homma}, we have that $p$ divides
$\epsilon_5$. More precisely we have the following. 
    \begin{lemma}\label{lemma3.2}
Either $\epsilon_5$ is a power of $p$ or $p=3$ and $\epsilon_5=6$.
    \end{lemma}
    \begin{proof} We can assume $\epsilon_5>5$. If $\epsilon_5$ is not a
power of $p$, by the $p$-adic criterion \cite[Corollary 1.9]{sv} we have 
$p\le 3$ and $\epsilon=6$.
    \end{proof} 
From Proposition \ref{prop3.1}, the case $\nu_4=\epsilon_5=6$ provides the
following bound:
  \begin{lemma}\label{lemma3.3}
Let $\cK$ be a complete arc of size $k$ such that $(2q+3)/4<k\le
m_2'(q)$. Suppose that $\cK$ admits an irreducible envelope such that 
$\nu_4=\epsilon_5=6$. Then $p=3$ and 
$$
k\le \frac{52}{53}q+\frac{44}{53}\, .
$$
  \end{lemma}
As in the previous case, for $q$ an even power of 3 and $k=m_2'(q)$
the case $\nu_4=\epsilon_5=6$ occur only for $q$ small. More precisely, we
have the following.  
  \begin{corollary}\label{cor3.2} Let $\cK$ be an arc of size $m_2'(q)$. 
Suppose that $q$ is an even power of $p$ and that $\cK$ admits an
irreducible envelope with $\nu_4=\epsilon_5=6$. Then $p=3$ and $q\le 3^6$.
  \end{corollary}
   \begin{proof} From the $p$-adic criterion \cite[Corollary 1.9]{sv},
$p=3$. Then from Proposition \ref{prop3.1} and (\ref{eq1.2}) we have
$$
q-\sq+1\le m_2'(q)\le 52q/53+44/53\, ,
$$
and the result follows.
  \end{proof} 
From now on we assume 
$$
\nu_4=\epsilon_5=\text{a power of $p$}\, . 
$$
Then, the bound 
   \begin{equation}\label{eq3.1}
k\le q-\frac{1}{4}\nu_4+\frac{7}{4}
   \end{equation}
in Proposition \ref{prop3.1} and
Segre's bound (\ref{eq1.1}) provide motivation to consider three cases
according as $\nu_4>\sq$, $\nu_4< \sq$, or $\nu_4=\sq$. 

{\bf 3.2.1. $\nu_4>\sq$.} 

Since $\nu_4$ is a power of $p$, then we
have that $\nu^2 \ge pq$ and so from (\ref{eq3.1}) the following
holds: 
   \begin{lemma}\label{lemma3.4} Let $\cK$ be a complete arc of size $k$
such that $(2q+4)/3<k\le m_2'(q)$. Suppose that $\cK$ admits an
irreducible envelope such that $\nu_4$ is a power of $p$ and
that $\nu_4>\sq$. Then
$$
k\le \begin{cases}
q-\frac{1}{4}\sqrt{pq}+\frac{7}{4} & \text{if $q$ is not a square}\, ,\\
q-\frac{1}{4}p\sq+\frac{7}{4} & \text{otherwise}\, .
     \end{cases}
$$
   \end{lemma}
If $q$ is a square and $k=m_2'(q)$, then $\nu_4>\sq$ can only occur in
characteristic 3:
   \begin{corollary}\label{cor3.3} Let $\cK$ be an arc of size $m_2'(q)$.
Suppose that $q$ is an even power of $p$ and that $\cK$ admits an
irreducible envelope with $\nu_4$ a power of $p$ and $\nu_4>\sq$.
Then $p=3$,
$\nu_4=3\sq$, and 
$$
k\le q-\frac{3}{4}\sq+\frac{7}{4}\, .
$$
   \end{corollary}
    \begin{proof} From Lemma \ref{lemma3.4} and $m_2'(q)\ge
q-\sq+1$ follow that $\sq(p-4)\le 3$ and so that $p=3$. From $\nu_4\le
2d-1$ and $2d\le 4t=4(q+2-m_2'(q))\le 4\sq+4$ we have that $\nu_4\le
4\sq+3$ and it follows the assertion on $\nu_4$. The bound on $k$ follows
from Lemma \ref{lemma3.4}.
   \end{proof}
{\bf 3.2.2. $\nu_4<\sq$.} 

Let 
$$
F(x):= (2x+32-q)/(4x+29)\, .
$$
Then the bound 
$$
k\le \frac{28+4\nu_4}{29+4\nu_4}q+\frac{32+2\nu_4}{29+4\nu_4}
$$ 
in Proposition \ref{prop3.1} can be written as 
   \begin{equation}\label{eq3.2}
k\le q+F(\nu_4)\, .
   \end{equation} 
For $x>0$, $F(x)$ is an increasing function so that 
$$
F(\nu_4)\le\begin{cases} 
F(\sqrt{q/p})= -\frac{1}{4}\sqrt{pq}+\frac{29}{16}p+\frac{1}{2}+R & 
\text{if $q$ is not a square}\, ,\\
F(\sq/p)=-\frac{1}{4}p\sq+\frac{29}{16}p^2+\frac{1}{2}+R &
\text{otherwise}\, ,
\end{cases}
$$
where
$$
R=\begin{cases} 
-\frac{841p-280}{16(4\sqrt{q/p}+29)} & \text{if $q$ is not a square}\, ,\\
-\frac{841p^2-280}{16(4\sq/p+29)}    & \text{otherwise}\, .
\end{cases}
$$
Then from (\ref{eq3.2}) and since $R<0$ we have the following.
   \begin{lemma}\label{lemma3.5} Let $\cK$ be a complete arc of size $k$
such that $(2q+3)/4<k\le m_2'(q)$. Suppose that $\cK$ admits an
irreducible envelope such that $\nu_4$ is a power of $p$ and
that $\nu_4<\sq$. Then
$$
k<\begin{cases}
q -\frac{1}{4}\sqrt{pq}+\frac{29}{16}p+\frac{1}{2} &  
\text{if $q$ is not a square}\, ,\\
q-\frac{1}{4}p\sq+\frac{29}{16}p^2+\frac{1}{2} &
\text{otherwise}\, .
\end{cases}
$$
    \end{lemma}
     \begin{corollary}\label{cor3.4} Let $\cK$ be a complete arc of size
$m_2'(q)$. Suppose that $q$ is an even power of $p$ and that $\cK$ admits
an irreducible envelope with $\nu_4$ a power of $p$ and $\nu_4<\sq$.
Then one of the following statements holds:
   \begin{enumerate}
\item $p=3$, $\nu_4=\sq/3$, and $m_2'(q)$ satisfies Lemma \ref{lemma3.5}.
\item $p=5$, $q=5^4$, $\nu_4=5$, and $m_2'(5^4)\le 613$;
\item $p=5$, $q=5^6$, $\nu_4=5^2$, and $m_2'(5^6)\le 15504$;
\item $p=7$, $q=7^4$, $\nu_4=7$, and $m_2'(7^4)\le 2359$.
   \end{enumerate}
     \end{corollary}
     \begin{proof} Let $q=p^{2e}$; so $e\ge 2$ as $p\le \nu_4<p^e$. 
From (\ref{eq1.2}) and Lemma \ref{lemma3.5} we have
that 
$$
(p-4)p^e/4<29p^2/16-0.5\, ,
$$
so that $p\in \{3,5,7,11\}$. 

Let $p=3$. If $\nu_4\le \sq/9$ (so $e\ge 4$), then from (\ref{eq1.2})
and $m_2'(q)\le q+F(\sq/9)$ we would have that 
$$
q-\sq+1\le q-9\sq/4+2357/16-67841/16(43^{e-2}+29)\, ,
$$
which is a contradiction for $e\ge 4$. 

Let $p=11$. Then $p^e\le 125$ and $e=2$ and $\nu_4=11$. Thus from
Proposition \ref{prop3.1} we have $m_2'(11^4)\le 11^4+F(11)$, i.e. 
$m_2'(11^4)\le 14441$. This is a contradiction since by (\ref{eq1.2}) we
must have $m_2'(11^4)\ge 14521$. This eliminates the possibility $p=11$. 

The other cases can be handled in an analogous way.
     \end{proof}
{\bf 3.2.3.  $\nu_4=\sq$.} 

In this case, according to (\ref{eq3.1}),
we just obtain Segre's bound (\ref{eq1.1}).
   \section{Irreducible envelopes of large complete arcs}\label{4}
Throughout this section we keep the notations of the previous section.
Here we study geometrical properties of irreducible envelopes associated
to large complete arcs in $\P^2(\fq)$, $q$ odd. To do so we use the bounds
obtained in \S3 and divide our study in two cases according as $q$ is a
square or not.

{\bf 1.  $q$ square.} 

Let $\cX$ be an irreducible envelope
associated to an arc of size $m_2'(q)$. Then from Lemma \ref{lemma2.3}, 
and Corollaries \ref{cor3.1}, \ref{cor3.2}, \ref{cor3.3}, \ref{cor3.4}, we
have the following
   \begin{proposition}\label{prop4.1} If $q$ is an odd square and
$q>43^2$, then 
$\cX$ is $\Sigma_1$-classical. The $\Sigma_2$-orders 
are $0,1,2,3,4, \epsilon_5$ and the 
$\fq$-Frobenius $\Sigma_2$-orders are $0,1,2,3,\nu_4$, with
$\epsilon_5=\nu_4$, where also one of the following holds:
   \begin{enumerate}
\item $\nu_4\in \{\sq/3,3\sq\}$ for $p=3$;
\item $(\nu_4,q)\in \{(5,5^4), (5^2,5^6), (7,7^4)\}$;
\item $\nu_4=\sq$ for $p\ge 5$.
   \end{enumerate}
   \end{proposition}
{\bf 2. $q$ non-square.} In this case there is no analogue to 
bound (\ref{eq1.2}). From Lemmas \ref{lemma3.1}, \ref{lemma3.3},
\ref{lemma3.4}, \ref{lemma3.5} and taking into consideration
({\ref{eq3.2}) we have the following.
   \begin{proposition}\label{prop4.2} Let $q>43^2$ and $q=p^{2e+1}$,
$e\ge 1$. Then, apart from the
values on $\nu_4$, 
the curve $\cX$, $\nu_4$ and $\epsilon_5$ are as in Proposition
\ref{prop4.1}. In this case 
$$
m_2'(q)>q-3\sqrt{pq}/4+7/4
$$
implies
  \begin{enumerate}
\item $\nu_4=\sqrt{q/p};$
\item $m_2'(q)<q-\sqrt{pq}/4+29p/16+1/2$.
  \end{enumerate}
  \end{proposition}
In particular our approach just
gives a proof of Segre's bound (\ref{eq1.1}) and Voloch's bound \cite{v2}.
However, both
propositions above show the type of curves associated to
large complete arcs. The study of such curves, for $q$ square and large enough,
allowed Hirschfeld and Korchm\'aros
\cite{hk1}, \cite{hk2} to improve Segre's bound (\ref{eq1.1}) to the bound in
(\ref{eq1.3}). For the sake of completeness we stress here the main ideas
from \cite{hk2} necessary to deal with Problem \ref{prob1.1}. Due to 
Proposition \ref{prop2.2}, 
the main strategy is to bound from above the 
number $2M_q+M_q'$ (which is defined via 
(\ref{eq2.1})). For instance, if one could prove that
  \begin{equation}\label{eq4.1}
2M_q+M_q'\le d(q-\sq+1)\, ,
   \end{equation}
where $d$ is the degree of the irreducible envelope whose normalization is
$\cX$, then from Proposition \ref{prop2.2} would follow immediately an
affirmative answer to Problem \ref{prob1.1}. However, since we know 
the answer to be negative for $q=9$ and $d\le 2t=2(q+2-m_2'(q))$, then one
can assume that $d$ is bounded by a linear function on $\sq$ and should  
expect to prove (\ref{eq4.1}) only under certain conditions on $q$. 
   \begin{lemma}\label{lemma4.1}
Let $q$ be an odd square. If (\ref{eq4.1}) holds true for $d\le
2\sq-\alpha$ with $\alpha\ge 0$, then $m_2'(q)<q-\sq+2+\alpha/2 $. In
particular, if (\ref{eq4.1}) holds 
true for $d\le 2\sq$, then the answer to Problem \ref{prob1.1} is
positive; that is,
$m_2'(q)=q -\sq+1$.
    \end{lemma}
    \begin{proof} If $m_2'(q)\ge q-\sq+2+\alpha/2$, then from $d\ge
2(q+2-m_2'(q)$ we would have that $d\le 2\sq-\alpha$ and so, from
Proposition \ref{prop2.2} and (\ref{eq4.1}), that $m_2'(q)\le q-\sq+1$, a
contradiction.
     \end{proof}
Now, in \cite{hk1}, Lemma \ref{lemma4.1} is proved  for $\alpha\ge \sq+3$,
i.e. whenever $d\le \sq-3$, and so (\ref{eq1.3}) follows. Recently,
Aguglia and Korchm\'aros \cite{a-k} proved a weaker version of
(\ref{eq4.1}) for $d=\sq-2$ and $q$ large enough, namely
$2M_q+M_q'<d(q-\sq/2-2)$. From this
inequality and Proposition \ref{prop2.2} one 
slightly improves (\ref{eq1.3}) to $m_2'(q)\le q-\sq/2-5/2$ whenever
$d=\sq-2$ and $q$ is large enough. Therefore the paper \cite{a-k}, as
well as \cite{hk1} 
or \cite{hk2}, is a good guide toward the proof of (\ref{eq4.1}) for
$\sq-2\le d\le 2\sq$.
\medskip

\begin{center}
{\bf APPENDIX: Background on Weierstrass points and Frobenius orders}
\end{center}
\medskip

In this section we summarize relevant material from St\"ohr-Voloch's paper
\cite{sv} concerning Weierstrass points and Frobenius orders.

Let $\cX$ be a projective geometrically irreducible non-singular algebraic
curve defined over $\bar\fq$ equipped with the action of the Frobenius
morphism $\frx$ over $\fq$. Let $p:={\rm char}(\fq)$. Let $\cD$ be a
base-point-free linear series $g^r_d$ 
on $\cX$ and assume that it is defined over $\fq$. Let $\pi:\cX\to
\P^r(\bar\fq)$ be the $\fq$-morphism associated to $\cD$. Then by
considering
the pullback of hyperplanes in $\P^r(\bar\fq)$ (via $\pi$)  one can
define, for each $P\in \cX$, a sequence of numbers
$j_0(P)=0<\ldots<j_r(P)$, called the {\em $(\cD,P)$-order sequence}. It
turns out that this sequence is the same, say
$\epsilon_0<\ldots<\epsilon_r$, for all but a finitely many points. This
constant sequence is called the {\em order sequence of $\cD$}. There
exists a divisor $R=R^\cD$, the so called {\em ramification divisor of
$\cD$}, such that the $\supp(R)$ is the set of points whose
$(\cD,P)$-orders are different from $(\epsilon_0,\ldots,\epsilon_r)$. The
curve is called {\em $\cD$-classical} with if 
$\epsilon_i=i$
for each $i$. The
following are the main properties of these invariants.
\begin{enumerate}
\item[(W1)] $j_i(P)\ge \epsilon_i$ for each $P$ and each $i$;
\item[(W2)] $v_P(R)\ge \sum_i(j_i(P)-\epsilon_i)$; equality holds iff
$\det(\binom{j_i(P)}{\epsilon_j})\not\equiv 0\pmod{p}$;
\item[(W3)] $\deg(R)=(2g-2)\sum_i\epsilon_i + (r+1)d$.
\end{enumerate}
Now to count $\fq$-rational points one looks for those points $P$ such
that $\pi(\frx(P))$ belongs to the osculating hyperplane at $P$. This led
to
the construction of a divisor $S=S^{\cD,q}$, the so
called {\em $\fq$-Frobenius divisor associated to $\cD$}, such that
\begin{enumerate}
\item[(F1)] $\cX(\fq)\subseteq \supp(S)$;
\item[(F2)] $\deg(S)=(2g-2)\sum_{i=0}^{r-1}\nu_i+(q+r)d$, where $\nu_0=0$
and $(\nu_1,\ldots,\nu_{r-1})$, called the {\em $\fq$-Frobenius orders of
$\cD$}, is a subsequence of $(\epsilon_1,\ldots,\epsilon_r)$.
\end{enumerate}
The curve is called {\em $\fq$-Frobenius classical with respect to $\cD$}
if $\nu_i=i$ for each $i$. In addition, for each $P\in \cX(\fq)$ holds
\begin{enumerate}
\item[(F3)] $\nu_i\le j_{i+1}(P)-j_1(P)$, $i=0,\ldots,r-1$;
\item[(F4)] $v_P(S)\ge \sum_{i=0}^{r-1}(j_{i+1}(P)-\nu_i)$.
\end{enumerate}
Hirschfeld and Korchm\'aros \cite{hk2} noticed that (F3) and (F4) even
holds for points in the set
$$
\cX_1(\fq):=\{P\in \cX: \pi(P)\in \pi(\cX)(\fq)\}\, .
$$
Therefore from (F3) and (F4) we have $v_P(S)\ge rj_1(P)$ for each $P\in
\cX_1(\fq)$ and hence we obtain the main result in \cite{sv}:
$$
\deg(S)/r\ge \#\cX_1(\fq)\ge \#\cX(\fq)\, .
$$ 
\smallskip

{\bf Acknowledgments.} The authors wish to thank James W.P. Hirschfeld and
G. Korchm\'aros for useful comments. This research was carried out with
the support of the Italian Ministry for Research and Technology (project
40\% ``Strutture geometriche, combinatorie e loro applicazioni"). Part of
this paper was written while Torres was visiting ICTP, Trieste-Italy
(June-July 1998) supported by IMPA/Cnpq-Brazil and ICTP.

\end{document}